\documentclass[a4paper,10pt]{article}

\usepackage[utf8]{inputenc}
\usepackage{amsmath}
\usepackage{amssymb}
\usepackage [hmargin={3cm,3cm},vmargin={2.4cm,2.6cm}]{geometry}
\usepackage{graphicx}

\newtheorem{thm}{Theorem}[section]

\newtheorem{Question}{Question}

\newtheorem{Conjecture}{Conjecture}

\newtheorem{prop}{Proposition}[section]

\newtheorem{lem}{Lemma}[section]

\newtheorem{rem}{Remark}[section]

\newtheorem{cor}{Corollary}[section]

\newtheorem{defin}{Definition}[section]

\newcommand{\R}{\mathbb{R}}

\newcommand{\pr}{\textit{Proof: \vskip5mm}}

\newcommand{\cqfd}{\begin{flushright}                  
			 $\Box$
                 \end{flushright}}

\title{A perturbation result for the Riesz transform}

\author{Baptiste Devyver\footnote{Laboratoire de Math\'{e}matiques Jean Leray, Universit\'{e} de Nantes (FRANCE); email: baptiste.devyver@univ-nantes.fr}}

\begin{document}

\maketitle


\begin{abstract}

We show a perturbation result for the boundedness of the Riesz transform: if $M$ and $M_0$ are complete Riemannian manifolds which are isometric outside a compact set, we give sufficient conditions so that the boundedness on $L^p$ of the Riesz transform on $M_0$ implies the boundedness on $L^p$ of the Riesz transform on $M$. 
\end{abstract}

\section{Introduction}

Let $(M,g)$ be a Riemannian manifold. The Riesz transform problem, namely giving conditions on $p$ and on the manifold such that the operator $d\Delta^{-1/2}$ -- the so-called Riesz transform -- is bounded on $L^p$, has recently undergone certain progresses. A pioneering result which goes back to 1985 is a theorem of D. Bakry \cite{Bakry} which asserts that if the Ricci curvature of $M$ is non-negative, then the Riesz transform on $M$ is bounded on $L^p$ for every $1<p<\infty$. However, it is only recently that some progresses have been made to understand the behaviour of the Riesz transform if some amount of negative Ricci curvature is allowed. A general question is the following:

\begin{Question}\label{Bakry}
What is the analogue of Bakry's result for manifolds with non-negative Ricci curvature (or with a ``small" amount of negative Ricci curvature) outside a compact set ?
\end{Question}
Here, the ``smallness" should be understood in an integral sense, i.e. $Ric_-\in L^r(\mathrm{d}\mu)$, for some value of $r$ and some measure $\mathrm{d}\mu$. A very partial answer has been provided in \cite{Coulhon-Zhang}, where it is shown essentially that if the Ricci curvature is smaller in an integral sense than a constant $\varepsilon$ (depending on the geometry of $M$), then the Riesz transform is bounded on $L^p$ for every $1<p<\infty$. However, this result is not entirely satisfying, since it does not say what happens if the integral of the Ricci curvature is bigger than the threshold $\varepsilon$: it does not thus cover the case of manifolds having non-negative Ricci curvature outside a compact set. Unlike manifolds with non-negative Ricci curvature, manifolds with non-negative Ricci curvature outside a compact set can exhibit several ends, as well as more complicated topology (although it is far from being clear how to quantify this), and it has been known for already some time that Bakry's result stated as such cannot hold for manifolds with non-negative Ricci curvature outside a compact set. Indeed, as is shown in \cite{Coulhon-Duong2}, the Riesz transform on $\R^n\#\R^n$, the connected sum of two Euclidean spaces, is not bounded on $L^p$ for $p>n$, although $\R^n\#\R^n$ has Riemannian curvature which vanishes outside a compact set. It was proved later in \cite{Carron-Coulhon-Hassell} that the Riesz transform on $\R^n\#\R^n$ is bounded on $L^p$ if and only if $1<p<n$ ($1<p\leq2$ if $n=2$).  In the same paper, the authors also prove that if the manifold has \textit{only one end} and is isometric outside the compact set to $\R^n$, then the Riesz transform is bounded on $L^p$ for all $1<p<\infty$. This was pushed further by C. Guillarmou and A. Hassell in \cite{Guillarmou-Hassell1} to study the Riesz transform on asymptotically conical manifolds: (a particular case of) their result is that when we glue together several conical manifolds, then if there are more than one end, the Riesz transform is bounded on $L^p$ iff $1<p<n$; and if there is only one end and the manifold is isometric outside a compact set to a conical manifold $M_0$, then the range of boundedness of the Riesz transform is the same as it is on $M_0$. The results cited above are in fact perturbation results for the Riesz transform, and let us reformulate them in the following way:

\begin{thm}\cite{Carron-Coulhon-Hassell},\cite{Guillarmou-Hassell1}\label{perturbation Carron-Coulhon-Hassell}
In the class of connected asymptotically Euclidean (or more generally asymptotically conical) manifolds of dimension $d$, the boundedness of the Riesz transform on $L^p$ is stable:

\begin{enumerate}
 \item Under ``gluing" (that is, connected sum construction), and change of both the metric and the topology on a compact set, if $1<p<d$.

 \item Under change of both the metric and the topology on a compact set, if $p\geq d$.
\end{enumerate}

\end{thm}
It is however a result very specific to the class of manifolds under consideration: the proofs rely on a precise study of the kernel of $d\Delta^{-1/2}$, using the difficult techniques of \textit{b-calculus}, for which we need a very precise description of the structure at infinity of both the manifold and the metric. There is thus no hope to generalize these proofs to general manifolds with non-negative Ricci curvature outside a compact set. Then G. Carron proved  in \cite{Carron4} a key perturbation result, which is more general. Let us say that on $(M,g)$ a \textit{Sobolev inequality} (of dimension $n>2$) holds if

\begin{equation}\label{Sobolev}\tag{$S_n$}
||f||_{\frac{2n}{n-2}}\leq C||\nabla f||_2,\,\forall f\in C_0^\infty(M)
\end{equation}
Then we define:
\begin{defin}
{\em 
The \textit{Sobolev dimension}  $d_S(M)$ is defined as the supremum of the set of $n$ such that the Sobolev inequality of dimension $n$ \eqref{Sobolev} is satisfied on $M$ (in the case where the Sobolev inequality is not satisfied for any $n$, we let by convention $d_S=-\infty$).
}
\end{defin}
The Sobolev dimension needs not be equal to the topological dimension of $M$, in fact one has only the inequality

$$d_S\geq \hbox{dim}(M)$$
(see \cite{Saloff-Coste1}). For asymptotically conical manifolds, the Sobolev dimension and the topological dimension coincide, but  $\mathbb{H}^n$, the hyperbolic space of dimension $n$, satisfies $d_S(\mathbb{H}^n)=+\infty$. Carron's perturbation result states as follows:

\begin{thm}\label{resultat Carron}\cite{Carron4}\\
Let $M_0$ and $M$ be a complete Riemannian manifolds (not necessarily connected), isometric outside a compact set, which satisfy $d_S>3$ and with Ricci curvature bounded from below. Assume that the Riesz transform on $M$ is bounded on $L^p$ for some $p\in (\frac{d_S}{d_S-1},d_S)$. Let $M$ be a complete Riemannian manifold isometric to $M_0$ at infinity, then the Riesz transform on $M$ is bounded on $L^p$.

\end{thm}
Here, isometric at infinity means that we can find two compact sets $K_0$ and $K$ (resp. of $M_0$ and of $M$), such that $M\setminus K$ is isometric to $M_0\setminus K_0$, and the fact that the manifolds are not supposed to be connected is to allow connected-sum constructions. This result extends $(1)$ of Theorem \ref{perturbation Carron-Coulhon-Hassell} to a much more general class of manifolds, namely to manifolds with Ricci curvature bounded from below, and satisfying a Sobolev inequality -- the dimension parameter up to which we can glue together two such manifolds while preserving the boundedness of the Riesz transform being the Sobolev dimension $d_S$. The fact that $d_S(M_0)=d_S(M)$ follows from \cite{Carron3}, Proposition 2.7. We see thus that, rather than the topological dimension, an important quantity from the point of view of the perturbation theory for the Riesz transform is the Sobolev dimension. \\

A way to rephrase Carron's result is that for $p<d_S$, the boundedness of the Riesz transform on $L^p$ is preserved under gluing and perturbation of the metric and the topology on a compact set. Thus, for example, the boundedness of the Riesz transform on $L^p$ for \textit{any} $1<p<\infty$ is preserved under gluing, perturbation of the topology and of the metric in the class of manifolds whose ends are isometric to $\mathbb{H}^n$ at infinity. However, when $d_S<\infty$, Carron's result does not say anything concerning the generalization of $(2)$ of Theorem \ref{perturbation Carron-Coulhon-Hassell}: explicitely, when $p\geq d_S$, what happens for the boundedness of the Riesz transform on $L^p$ if we start with a manifold with one end, and we change both the metric and the topology on a compact set, without making any gluing, i.e. preserving the fact that the manifold has only one end? \\

Let us mention at this point a perturbation result of Coulhon and Dungey \cite{Coulhon-Dungey} which investigates what happens for the Riesz transform if we change the metric and the Riemannian measure. Under quite mild conditions on the perturbation, they show that the boundedness on $L^p$ of the Riesz transform is preserved under a change of metric and of measure, \textit{for any} $1<p<\infty$. However, their main assumption is that the underlying manifold is the same, that is they allow no change of topology \textit{at all}, and their method relies crucially on this assumption. As a consequence, it is not possible from their result to get either $(1)$ or $(2)$ of Theorem \ref{perturbation Carron-Coulhon-Hassell}, even for the case of the Euclidean space.\\\\
In the paper \cite{Devyver-Riesz1}, we use Carron's perturbation result to answer Question \ref{Bakry} for the case $p<d_S$: under the assumptions that $d_S>3$, that the negative part of the Ricci curvature is in $L^{\frac{d_S}{2}-\varepsilon}\cap L^\infty$, and that the volume of balls of big radius $R$ is comparable to $R^{d_S}$, we show that the Riesz transform is bounded on $L^p$ for $1<p<d_S$. If in addition there are no non-zero $L^2$ harmonic $1-$forms, we also prove that the Riesz transform is bounded on $L^p$ for all $1<p<\infty$. However, it is expected that this last assumption is too strong to get the boundedness on the whole $(1,\infty)$, more precisely in \cite{Devyver-Riesz1} we made the following conjecture:

\begin{Conjecture}\label{conj}
If $M$ is a manifold satisfying a Sobolev inequality, having non-negative Ricci curvature outside a compact set and only one end, then the Riesz transform on $M$ is bounded on $L^p$ for every $1<p<\infty$.
 \end{Conjecture}
In other words, is the presence of several ends the only obstruction in this class of manifolds to the boundedness of the Riesz transform on $L^p$ for all $1<p<\infty$? Motivated by this conjecture, we generalize in this article both Theorem \ref{perturbation Carron-Coulhon-Hassell} and Theorem \ref{resultat Carron}. We will assume that the manifold satisfies a Sobolev inequality so that $d_S$, the Sobolev dimension, is greater than $2$, and we will be interested in extending the mentionned perturbation results Theorems \ref{perturbation Carron-Coulhon-Hassell} and \ref{resultat Carron} to the case where $p\geq d_S$. First, we define the hyperbolic dimension of $M$ to be (see section 1)

\begin{defin}

{\em 
The \textit{hyperbolic dimension} $d_H(M)$ of $M$ is the supremum of the set of $p$ such that $M$ is $p-$hyperbolic.

}

\end{defin}
Our main result shows first that $d_H$ -- and not $d_S$ as Carron's result seems to indicate -- is the relevant quantity to be considered from the point of view of the behaviour of the Riesz transform under gluing; and secondly, we are able to generalize (2) of Theorem \ref{perturbation Carron-Coulhon-Hassell} under much more general assumptions. Our result writes:

\begin{thm}\label{Riesz p-hyperbolique}

Let $M$, $M_0$ be two Riemannian manifolds (not necessarily connected), isometric outside a compact set, whose Ricci curvature is bounded from below and satisfying $d_S>2$. We assume that the Riesz transform on $M_0$ is bounded on $L^p$ for $p\in [p_0,p_1)$ with $\frac{d_S}{d_S-1}<p_0\leq 2$ and $p_1>\frac{d_S}{d_S-2}$. Then the Riesz transform on $M$ is bounded on $L^p$ for $p\in\left[p_0,\min(d_H(M),p_1)\right)$. If furthermore $M$ has only one end, then the Riesz transform on $M$ is bounded on $L^p$ for $p\in\left[p_0,p_1\right)$.

\end{thm}
We now make a certain number of comments about this result:
\begin{rem}
{\em 
\begin{enumerate}
\item We will prove in section 1 (Proposition \ref{dimension hyperbolic}) that if the Riesz transform on $M$ is bounded on $L^p$ for $p\in \left(\frac{d_S}{d_S-1},2\right]$, then

$$d_S\leq d_H,$$
so that under this mild assumption our result indeed generalizes Carron's result (up to endpoints of the range of boundedness). Our result says that $d_H$, and not $d_S$, is the relevant quantity to be considered when we perform a gluing. However, due to the fact that the behaviour of the Riesz transform is not known for many examples, we do not know (although we think there exists) an example of a manifold $M_0$  for which $p_1>d_S$ and $d_H>d_S$. Nonetheless, we will see in Corollary \ref{Riesz non borne somme connexe p-hyperbolique} an application using $d_H$ and not $d_S$.

\item In the case where $M$ has only one end, this result extends point $(2)$ of Theorem \ref{perturbation Carron-Coulhon-Hassell} to the class of manifolds satisfying a Sobolev inequality. This provides evidence in favour of Conjecture \ref{conj}, and it could be also a necessary tool to prove it, in the same way that we used Carron's result \cite{Carron4} in \cite{Devyver-Riesz1} in order to prove boundedness of the Riesz transform on $L^p$ for $p<d_S$.

\item We expect that the hypothesis that $M_0$ satisfies a Sobolev inequality is too strong. A more reasonnable hypothesis would be that $M$ satisfies the \textit{relative Faber-Krahn inequality}, which is equivalent (see \cite{Grigor'yan5}) to the fact that the volume form on $M$ is doubling and that the heat kernel of the Laplacian satisfies a Gaussian upper estimate.

\item We are not able to treat the upper endpoint of the range of boundedness. However, in all the known cases, the range of boundedness of the Riesz transform is an \textit{open} interval, i.e. the Riesz transform is not bounded at the endpoints, and thus our limitation is not so disturbing. We also need to assume the the technical condition $p_1>\frac{d_S}{d_S-2}$, which is satisfied in most of the cases (and in all the interesting cases covered by Carron's result, when $d_S\geq4$).

\item Recall that in Carron's result, one needs to assume $d_S>3$. In our result, we can allow $d_S=3$, but in this case

$$\frac{d_S}{d_S-2}=d_S,$$ 
therefore we cannot deduce from our result that the Riesz transform on the connected sum of two copies of $R^3$ is bounded on $L^p$ for $1<p<3$ (which is true and was proved in \cite{Carron-Coulhon-Hassell}).

\end{enumerate}
}
\end{rem}
Theorem \ref{Riesz p-hyperbolique} has a certain number of interesting corollaries, which we describe now. The first three of them follow from Theorem (\ref{Riesz p-hyperbolique}) with the hypothesis ``$M$ has only one end", and the last one uses the hyperbolic dimension $d_H$. First, we recover a particular case of a result of C. Guillarmou and A. Hassell \cite{Guillarmou-Hassell1} on asymptotically conical manifolds, without using the heavy machinery of \textit{b-calculus} -- as in \cite{Guillarmou-Hassell1}, this uses H.Q. Li's result \cite{Li1} about the Riesz transform on conical manifolds.
\begin{cor}\label{asymptotic conic}

Let $M$ be a complete Riemannian manifold, isometric outside a compact set to a conical manifold $M_0=\R_+^\star\times N$, with $(N,h)$ connected and compact of dimension $n-1$ -- that is, $M_0$ is endowed with the metric $g=dt^2+t^2h$. Let $\lambda_1$ be the first non-zero eigenvalue of the Laplacian on $N$, and let

$$p_0:=\frac{n}{\frac{n}{2}-\sqrt{\lambda_1+\left(\frac{n-2}{2}\right)^2}}$$
(with $p_0=\infty$ if $\lambda_1\geq n-1$ by convention). Then if $n\geq3$, the Riesz transform on $M$ is bounded on $L^p$ when $1<p<p_0$, and is unbounded on $L^p$ when $p>p_0$.

\end{cor}
Furthermore, we also have the following new results:

\begin{cor}\label{Riesz Ricci positive hors d'un compact}

Let $M$ be a complete Riemannian manifold with one end, isometric outside a compact set to a manifold with non-negative Ricci curvature. We assume that on $M$ we have the following volume estimate: there is $x_0\in M$ and $\nu>2$ such that

$$V(x_0,R)\geq C R^\nu,\,\forall R>0,$$
then the Riesz transform on $M$ is bounded on $L^p$ for all $\frac{\nu}{\nu-1}<p<\infty$.

\end{cor}

\begin{cor}\label{Lie group}

Let $M_0$ be a connected and simply connected nilpotent Lie group, and let $M$ be isometric outside a compact set to $M_0$. Then the Riesz transform on $M$ is bounded on $L^p$ for every $1<p<\infty$.

\end{cor}
Finally, we have the following corollary, which is also new:

\begin{cor}\label{Riesz non borne somme connexe p-hyperbolique}

Let $n\geq3$, and let $N$ be a manifold which is $q-$hyperbolic for some $q>n$, and which has Ricci curvature bounded from below. Then the Riesz transform on $M=N\#\R^n$, the connected sum of $N$ and $\R^n$, is not bounded on $L^p$ for $n<p<\infty$. In particular, the Riesz transform on the connected sum $\R^n\#\mathbb{H}^n$ of an Euclidean space and a hyperbolic space is not bounded on $L^p$ for $n<p<\infty$.

\end{cor}
The organization of this article is as follows: in section 1, we mostly review some classical results concerning the notion of $p-$hyperbolicity. In section 2, we prove Theorem \ref{Riesz p-hyperbolique} and its corollaries.

\section{About p-hyperbolicity}

In this section we recall some notions concerning p-hyperbolicity that we will need in the sequel. References for this section are \cite{Coulhon-Saloff-Coste-Holopainen} and \cite{Goldstein-Troyanov}. We will assume that the manifold is smooth, so that local elliptic theory applies. In particular, we will make use of the local Sobolev injections, of the trace theorems and of Poincar\'{e} inequalities for bounded domains. For references on this, see \cite{Saloff-Coste1} and \cite{Taylor}. Let us fix $1< p<\infty$.

\begin{defin}
{\em 
We say that a Riemannian manifold $(M,g)$ is \textit{p-hyperbolic} if for every non-empty, relatively compact open subset $U$ of $M$, there exists a constant $C_U$ such that

$$\int_U |f|^p\leq C_U\int_M|\nabla f|^p,\,f\in C_0^\infty(M).$$
}
\end{defin}
As in the case $p=2$, we have the following Proposition:

\begin{prop}\label{p-parabolicite}

$(M,g)$ is p-hyperbolic if and only if there exists \textbf{some} non-empty, relatively compact open subset  $U$ of $M$ and a constant $C_U$ such that

$$\int_U |f|^p\leq C_U\int_M|\nabla f|^p,\,f\in C_0^\infty(M).$$

\end{prop}
We write the proof for the reader's convenience.\\
\pr It is enough to show that for every smooth connected open set $W$ containing $U$, there exists $C_W$ such that

$$\int_W |f|^p\leq C_W\int_M|\nabla f|^p,\,f\in C_0^\infty(M).$$
We will need the following Lemma:

\begin{lem}\label{trou spectral p-laplacien}

For every relatively compact open sets $\Omega_1$, $\Omega_2$, such that $\Omega_1\subset\subset \Omega_2$ and such that $\Omega:=\Omega_2\setminus \overline{\Omega}_1$ is connected, there exists a constant  $C_\Omega$ such that

\begin{equation}\label{spectre p-laplacien}
||f||_p\leq C_\Omega||\nabla f||_p,\,\forall f\in C_{D-N}^\infty(\Omega),
\end{equation}
where $C_{D-N}^\infty(\Omega)$ is the set of smooth functions on $\Omega$ taking value $0$ on $\partial \Omega_1$ (the index $D-N$ stands for  ``Dirichlet-Neumann").

\end{lem}
Let us assume for a moment the result of the Lemma, and let us conclude the proof of Poposition (\ref{p-parabolicite}). Let $V$ be a non-empty, open set such that $V\subset\subset U$ and such that $W\setminus V$ is connected, and let $\rho$ be a smooth function whose support is included in $U$, such that $\rho\equiv1$ on $V$. Then

$$||f||_{L^p(W)}\leq ||\rho f||_{L^p(W)}+||(1-\rho)f||_{L^p(W)}.$$
Since $||\rho f||_{L^p(W)}=||\rho f||_{L^p(U)}$, we have by hypothesis

$$||\rho f||_{L^p(W)}\leq C_U||\nabla(\rho f)||_p\leq C_U\left(||f\nabla\rho||_p+||\rho\nabla f||_p\right).$$
On an other end, $||\rho\nabla f||_p\leq ||\rho||_\infty||\nabla f||_p$, and by hypothesis, since the support of $\nabla \rho$ is contained in $U$,

$$||f\nabla\rho||_p\leq||\nabla\rho||_\infty ||f||_{L^p(U)}\leq C||\nabla f||_p.$$
It remains to treat the term $||(1-\rho)f||_{L^p(W)}$. We apply Lemma (\ref{trou spectral p-laplacien}) with $\Omega=W\setminus V$, to obtain

$$||(1-\rho)f||_{L^p(W)}\leq C||\nabla\left((1-\rho)f\right)||_p\leq C\left(||\nabla(\rho f)||_p+||\rho||_\infty||\nabla f||_p\right),$$
and we bound as before $||\nabla(\rho f)||_p$ by $C||\nabla f||_p$.
\cqfd
\textit{Proof of Lemma (\ref{trou spectral p-laplacien}):} By contradiction, if there exists a sequence of functions $f_n\in C_{D-N}^\infty$ such that $||f_n||_{L^p}=1$, and $||\nabla f_n||_{L^p}\rightarrow 0$. Since $W^{1,p}(\Omega)$ is reflexive for $1<p<\infty$, up to the extraction of a subsequence we can assume that the sequence $(f_n)_n$ converges weakly to $f$ in $W^{1,p}(\Omega)$. But we have the compact Sobolev injection $W^{1,p}(\Omega)\hookrightarrow L^p$, therefore $(f_n)_n$ converges strongly in $L^p$, and as a consequence converges strongly to $f$ in $W^{1,p}(\Omega)$. The function $f$ then satisfies $\nabla f=0$ in the weak sense, and this implies that $\nabla f=0$ strongly, hence $f$ is constant since $\Omega$ is connected. In addition, the trace theorem for $W^{1,p}$ shows that $f|_{\partial \Omega_1}=0$, and therefore $f$ is zero. This contradicts the fact that $||f||_p=1$.
\cqfd
We will also use another caracterisation of $p-$hyperbolicity. Let us define first:

\begin{defin}
{\em 
If $U$ is a non-empty, relatively compact open subset of $M$, we define its $p-$\textit{capacity} by

$$\hbox{Cap}_p(U)=\inf\left\{\int_{M}|\nabla u|^p : u\in C_0^\infty\hbox{ such that }u|_U\geq 1\right\}=\inf\left\{\int_{M}|\nabla u|^p : u\in C_0^\infty\hbox{ such that }u|_U\equiv1\right\}.$$
}
\end{defin}
The last inequality in this definition follows from the fact that the ``truncation" of a function $u$ up to height $1$ on $U$ decreases the energy  $\int_M |\nabla u|^p$. For a detailed proof, see \cite{Goldstein-Troyanov}, Corollary 7.5. With this definition, we have the following caracterisation of the $p-$hyperbolicity:

\begin{thm}\label{hyperbolicity and capacity}

$M$ is $p-$hyperbolic if and only if the $p-$capacity of some (all) non-empty, relatively compact open set is non-zero.

\end{thm}
For a proof, see \cite{Troyanov}, Proposition 1.

\begin{rem}\label{volume hyperbolic}
{\em 
With the result of Theorem (\ref{hyperbolicity and capacity}), it is easy to see that if $M$ is $p-$hyperbolic for some $1<p<\infty$, then $M$ has infinite volume.
}
\end{rem}

\begin{cor}\label{somme connexe p-hyperbolique}

A Riemannian manifold $(M,g)$ is $p-$hyperbolic if and only if one of its ends is $p-$hyperbolic.

\end{cor}
\pr It is enough to find a non-empty, relatively compact open subset $\Omega$ of $M$, whose $p-$capacity is non-zero. We take $\Omega$ such that $M\setminus \Omega=M_1\setminus B_1\sqcup\ldots\sqcup M_k\setminus B_k$, the $M_i$ being the (closed) ends of $M$, and the $B_i$ being non-empty, relatively compact open subsets of $M_i$. Using the fact that the $p-$capacity of a non-empty, relatively compact open subset $U$ is equal to

$$\inf\left\{\int_{M\setminus U}|\nabla u|^p : u\in C_0^\infty\hbox{ such that }u|_U\equiv 1\right\},$$
we see that 

$$\hbox{Cap}_p(\Omega)=\sum_{i=1}^k\hbox{Cap}_p^{M_i}(B_i).$$
By hypothesis, one of the $M_i$ is $p-$hyperbolic ($M_1$ for example), which implies

$$\hbox{Cap}^{M_1}_p(B_1)>0,$$
and therefore

$$\hbox{Cap}_p(\Omega)>0.$$
\cqfd 
The main result of this section is the following link between $p-$hyperbolicity and Riesz transform:

\begin{prop}\label{laplacien non p-parabolique}

Let $M$ be a Riemannian manifold, which is $p-$hyperbolic for a certain $1<p<\infty$. We assume that the Riesz transform on $M$ is bounded on $L^p$. Then

$$\Delta^{-1/2} : L^p\rightarrow L^p_{loc},$$
is a bounded operator. Conversely, if the Riesz transform is bounded on $L^{q}$, $q$ being the dual exponent of $p$, and if 

$$\Delta^{-1/2} : L^p\rightarrow L^p_{loc},$$
is a bounded operator, then $M$ is $p-$hyperbolic.

\end{prop}
\pr Recall that the domain $L^p$ of $\Delta^{1/2}$ is defined as the set of functions $h$ in $L^p$ such that $\frac{e^{-t\sqrt{\Delta}}h-h}{t}$ has a limit in $L^p$ when $t$ tends to $0$. We will first prove the following Lemma:

\begin{lem}\label{densite p-hyperbolique}

For $1<p<\infty$, $C_0^\infty(M)$ is contained in the domain $L^p$ of $\Delta^{1/2}$, and $\Delta^{1/2}C_0^\infty$ is dense in $L^p$. Furthermore, if $u\in C_0^\infty$, then

$$\Delta^{-1/2}\Delta^{1/2}u=u.$$

\end{lem}
\textit{Proof of Lemma (\ref{densite p-hyperbolique}):} If $f\in C_0^\infty(M)$, we write

$$\Delta^{1/2}f=\Delta^{-1/2}\Delta f=\int_0^\infty e^{-t\sqrt{\Delta}}\Delta fdt,$$
and we separate the integral in $\int_0^1+\int_1^\infty=I_1+I_2$. In order to bound the $L^p$ norm of $I_1$, we use the fact that $\Delta f\in L^p$ and that $||e^{-t\sqrt{\Delta}}||_{p,p}\leq 1$, which yields

$$||I_1||_p\leq ||\Delta f||_p.$$
For $I_2$, we use the analyticity of $e^{-t\sqrt{\Delta}}$ on $L^p$, which implies that

$$\left|\left|\Delta e^{-t\sqrt{\Delta}}\right|\right|_{p,p}\leq \frac{C}{t^2}.$$
Consequently, we obtain

$$||I_2||_p\leq C||f||_p,$$
which gives that $\Delta^{1/2}f\in L^p$.\\

Let us now show that $\Delta^{1/2}C_0^\infty$ is dense in $L^p$. First, $(\Delta+1)C_0^\infty$ is dense in $L^p$: indeed, if $f\in L^q$ is orthogonal to $(\Delta+1)C_0^\infty$ (where $q$ is the conjugate exponent of $p$), then we have in the weak sense

$$(\Delta+1)f=0,$$
and this implies by a result of S.T. Yau (see Theorem 4.1 of \cite{Pigola-Rigoli-Setti}) that $f$ is constant, then that $f$ is zero since $M$ is of infinite volume by Remark (\ref{volume hyperbolic}). So $(\Delta+1)C_0^\infty$ is dense in $L^p$. Then,  $\Delta^{1/2}\left(\Delta+1\right)^{-1}$ is a bounded operator on $L^p$: to see this, we write

$$\Delta^{1/2}\left(\Delta+1\right)^{-1}=\int_0^\infty \Delta^{1/2}e^{-t\Delta}e^{-t}dt,$$
and we use the analyticity of $e^{-t\Delta}$ to say that

$$\left|\left|\Delta^{1/2}e^{-t\Delta}\right|\right|_{p,p}\leq \frac{C}{\sqrt{t}},\,\forall t>0.$$
Now, we write

$$\Delta^{1/2}C_0^\infty=\Delta^{1/2}\left(\Delta+1\right)^{-1}(\Delta+1)C_0^\infty,$$
and since $(\Delta+1)C_0^\infty$ is dense in $L^p$, and that $\Delta^{1/2}\left(\Delta+1\right)^{-1}$ is continuous on $L^p$, we have to see that the range of $\Delta^{1/2}\left(\Delta+1\right)^{-1}$ is dense in $L^p$. But $\left(\Delta+1\right)^{-1}L^p=\mathcal{D}_p(\Delta)$, the domain $L^p$ of the Laplacian. So we have to see that $\Delta^{1/2}\mathcal{D}_p(\Delta)$ is dense in $L^p$. But $\mathcal{D}_p(\Delta)$ contains $\Delta^{1/2}C_0^\infty$ by the first part of the Lemma: indeed, if $g\in C_0^\infty$, 

$$\Delta (\Delta^{1/2}g)=\Delta^{1/2}(\Delta g),$$
and this is in $L^p$ since $\Delta g\in L^p$. Therefore $\Delta^{1/2}\mathcal{D}_p(\Delta)$ contains

$$\Delta^{1/2}\Delta^{1/2}C_0^\infty=\Delta C_0^\infty,$$
which is dense in $L^p$ again by Yau's result.\\\\
It remains to show that when $u\in C_0^\infty$, then

$$\Delta^{-1/2}\Delta^{1/2}u=u.$$
We write

$$\begin{array}{rcl}
\Delta^{-1/2}\Delta^{1/2}u&=&\int_0^\infty e^{-t\sqrt{\Delta}}\Delta^{1/2} u\,dt\\\\
&=&\int_0^\infty-\frac{d}{dt}\left(e^{-t\sqrt{\Delta}}u\right)\,dt\\\\
&=&u-\lim_{t\to\infty}e^{-t\sqrt{\Delta}}u
\end{array}$$
By the spectral theorem, $\lim_{t\to\infty}e^{-t\sqrt{\Delta}}u$ converges to the projection of $u$ on the $L^2$ kernel of $\Delta$. But by Yau's above-mentionned result and the fact that $M$ has infinite volume, the $L^2$ kernel of $\Delta$ is reduced to $\{0\}$, and therefore

$$\Delta^{-1/2}\Delta^{1/2}u=u.$$
\cqfd
Now, we come back to the proof of Proposition (\ref{laplacien non p-parabolique}). We consider the first part of the Proposition. Let $\Omega$ be a non-empty, open, relatively compact set in $M$. The fact that the Riesz transform is bounded on $L^p$ is equivalent to the inequality

$$||\nabla u||_p\leq C||\Delta^{1/2}u||_p,\,\forall u\in C_0^\infty.$$
Since $M$ is $p-$hyperbolic, we also have the inequality

$$||u||_{L^p(\Omega)}\leq C||\nabla u||_p,\,\forall u\in C_0^\infty.$$
Combining these two inequalities, we obtain

$$||u||_{L^p(\Omega)}\leq C||\Delta^{1/2}u||_p,\,\forall u\in C_0^\infty.$$
Fix $u\in C_0^\infty$, and define $v=\Delta^{1/2}u$. Using the fact that since $u\in C_0^\infty$, by Lemma \ref{densite p-hyperbolique}

$$\Delta^{-1/2}\Delta^{1/2}u=u,$$
we have that $v$ is in the $L^p$ domain of $\Delta^{-1/2}$, and moreover $\Delta^{-1/2}v=u$. We thus obtain

$$||\Delta^{-1/2}v||_{L^p(\Omega)}\leq C||v||_p.$$
This is true for every $v\in \Delta^{1/2}C_0^\infty$, but by Lemma \ref{densite p-hyperbolique} $\Delta^{1/2}C_0^\infty$ is dense in $L^p$, and thus we obtain that 

$$||\Delta^{-1/2}v||_{L^p(\Omega)}\leq C||v||_p,\,\forall v\in L^p,$$
which is the result of the first part.\\

For the converse, we start with the assumption that there is a constant $C$ and a non-empty, open, relatively compact set $\Omega$ such that  for every $v\in L^p$,

$$||\Delta^{-1/2}v||_{L^p(\Omega)}\leq C ||u||_p.$$
Apply this to $v:=\Delta^{1/2}u$ for $u\in C_0^\infty(M)$ (which is licit by Lemma \ref{densite p-hyperbolique}), and using that $\Delta^{-1/2}v=u$ gives

$$||u||_{L^p(\Omega)}\leq C||\Delta^{1/2}u||_p,\,\forall u\in C_0^\infty(M).$$
But it is well-known that the boundedness of the Riesz transform on $L^q$ gives following the dual inequality: there is a constant $C$ such that

$$||\Delta^{1/2}u||_p\leq C||\nabla u||_p,\,\forall u\in C_0^\infty(M).$$
As a consequence, we get

$$||u||_{L^p(\Omega)}\leq C||\nabla u||_p,\,\forall u\in C_0^\infty(M),$$
i.e. $M$ is $p-$hyperbolic.

\cqfd
To conclude this section, we prove an inequality, announced in the introduction, involving the hyperbolic dimension and the Sobolev dimension. First, recall that the definition:

\begin{defin}
The \textit{hyperbolic dimension} $d_H$ of $M$ is defined as the supremum of the set of $p$ such that $M$ is $p-$hyperbolic. 
\end{defin}
Notice that (up to the author's knowledge) it is not known in full generality that the set of $p$ such that $M$ is $p-$hyperbolic is an interval; of course, by Proposition \ref{laplacien non p-parabolique}, this is true if the Riesz transform on $M$ is bounded on $L^p$ for $1<p<\infty$. Furthermore, by Corollary \ref{somme connexe p-hyperbolique}, if $M$ and $M_0$ are isometric outside a compact set, then 

$$d_H(M_0)=d_H(M).$$

\begin{prop}\label{dimension hyperbolic}
Let $M$ satisfying $d_S>2$, and assume that the Riesz transform on $M$ is bounded on $L^p$ for $p\in \left(\frac{d_S}{d_S-1},2\right]$. Then

$$d_H\geq d_S.$$
More precisely, $M$ is $p-$hyperbolic for every $2\leq p<d_S$.
\end{prop}
\pr Denote $d=d_S$, and let $1<p<d$. By Varopoulos \cite{Varopoulos},

$$\Delta^{-1/2} : L^p\rightarrow L^{\frac{dp}{d-p}},$$
and in particular 

$$\Delta^{-1/2} : L^p\rightarrow L^p_{loc}.$$
Let $\Omega$ be an non-empty, open, relatively compact set of $M$, then for every $u\in L^p$,

$$||\Delta^{-1/2}u||_{L^p(\Omega)}\leq C ||u||_p.$$
Applying this inequality to $u=\Delta^{1/2}v$ for $v\in C_0^\infty$ (this is licit since we have $\Delta^{1/2}C_0^\infty\subset L^p$ by Lemma \ref{densite p-hyperbolique}) yields

$$||v||_{L^p(\Omega)}\leq C||\Delta^{1/2}v||_p,\,\forall v\in C_0^\infty.$$
But it is well-knwon that the boundedness of the Riesz transform on $L^q$ implies the dual inequality

$$||\Delta^{1/2}v||_{q'}\leq C||\nabla v||_{q'},\,\forall v\in C_0^\infty,$$
where $q'$ is the dual exponent of $q$. By hypothesis, the Riesz transform is bounded on $L^{p'}$, so that, using the dual inequality for $q=p'$, we obtain

$$||v||_{L^p(\Omega)}\leq C||\nabla v||_p\,\forall v\in C_0^\infty,$$
which is exactly saying that $M$ is $p-$hyperbolic.

\cqfd

\section{Proof of the main results}

This section is devoted to the proof of Theorem \ref{Riesz p-hyperbolique} and its corollaries, announced in the introduction. We will extend the proof of Theorem (\ref{resultat Carron}) in \cite{Carron4}, to get rid of the condition $p<d_S$. For the convenience of the reader, we have divided the proof in several subsections. First, in subsection 1, we introduce several definitions and notations. In subsection 2, we recall the construction of \cite{Carron4}. In subsection 3, we prove Theorem \ref{Riesz p-hyperbolique} in the case of several ends. In subsection 4, we prove Theorem \ref{Riesz p-hyperbolique} in the case of one end. Finally, in subsection 5, we prove the corollaries of Theorem \ref{Riesz p-hyperbolique}.\\

\subsection{Definitions and notations}
\textit{Notation:} we will write $\mathfrak{s}(f)$ for the support of $f$.\\\\
Let $K_1$ be a compact set with smooth boundary such that $M\setminus K_1$ is isometric to the complement of a compact set of $M_0$, and $K_2$, $K_3$ compact sets with smooth boundaries such that $K_1\subset K_2\subset K_3$ and such that $K_i$ is contained in the interior of $K_j$ for $i<j$. \\

We define $\Omega:=M\setminus K_1$. Let $(\rho_0,\rho_1)$ be a partition of unity such that $\rho_1|_{K_1}\equiv1$,  $\mathfrak{s}(\rho_0)\subset\Omega$ and  $\mathfrak{s}(\rho_1)\subset K_2$. We also take $\varphi_0$ and $\varphi_1$ two smooth functions, such that $\mathfrak{s}(\varphi_0)\subset\Omega$, $\mathfrak{s}(\varphi_1)\subset K_3$ and such that $\varphi_i\rho_i=\rho_i$ for $i=1,2$. Furthermore, we assume that $\varphi_1 |_{K_2}\equiv 1$. \\

We will denote by $A$ the closure of a relatively compact, smooth open subset containing $\mathfrak{s}(d\varphi_0)$. We can arrange so that the distance between $A$ and  $\mathfrak{s}(\rho_0)$ is non-zero. Moreover, we can arrange so that $A$ is a disjoint union of connected ``annuli" $A_i$, each annulus corresponding to an end of $M_0$.
\begin{center}
\includegraphics[scale=1]{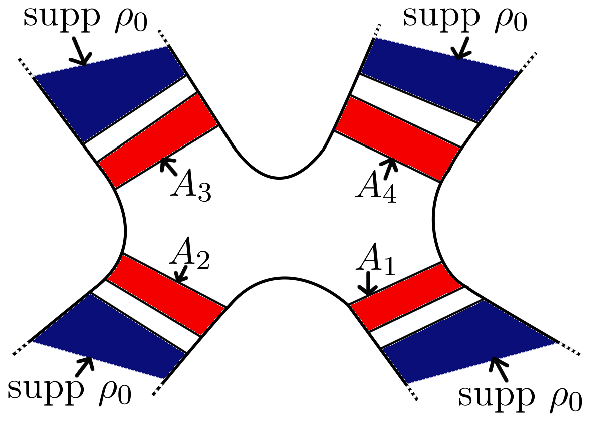}
\end{center}

\subsection{About Carron's proof of Theorem \ref{resultat Carron}}
The idea of G. Carron to prove Theorem \ref{resultat Carron} is to build a parametrix for $e^{-t\sqrt{\Delta}}$. Then by the formula 

$$\Delta^{-1/2}=\int_0^\infty e^{-t\sqrt{\Delta}}dt,$$
the parametrix for $e^{-t\sqrt{\Delta}}$ integrated in time yields a parametrix for $\Delta^{-1/2}$, and by differentiation of the Riesz transform $d\Delta^{-1/2}$. Therefore Carron's proof is in two steps: first, the construction of a good parametrix for $e^{-t\sqrt{\Delta}}$, such that when integrated and differentiated, it will yield a parametrix bounded on $L^p$ for the Riesz transform. And secondly, one needs to prove that the error term between the parametrix and the Riesz transform is also bounded on $L^p$.\\\\
Explicitely, Carron takes for the parametrix of $e^{-t\sqrt{\Delta}}$:

$$\mathfrak{E}(u)(s)=\varphi_0e^{-s\sqrt{\Delta_0}}\rho_0u+\varphi_1e^{-s\sqrt{\Delta_1}}\rho_1u,\,\forall u\in C_0^\infty(M),$$
where $\Delta_0$ is the Laplacian on $M_0$, and $\Delta_1$ is the Laplacian on $K_3$ with Dirichlet boundary conditions. Then we have the following formula:

$$e^{-t\sqrt{\Delta}}u=\mathfrak{E}(u,t)-G\left(-\frac{\partial^2}{\partial t^2}+\Delta\right)\mathfrak{E}(u),$$
where $G$ is the Green operator of $\left(-\frac{\partial^2}{\partial t^2}+\Delta\right)$ over $\R_+\times M$ with Dirichlet boundary conditions on the boundary $\{0\}\times M$. Indeed, at $t=0$, $e^{-t\sqrt{\Delta}}u=\mathfrak{E}(u,t)=u$. The term $G\left(-\frac{\partial^2}{\partial t^2}+\Delta\right)\mathfrak{E}(u)$ is the error term in the parametrix of $e^{-t\sqrt{\Delta}}$. When integrated and differentiated, the above parametrix for $e^{-t\sqrt{\Delta}}$ yields a parametrix for the Riesz transform, which is explicitely

$$\mathfrak{R}:=\sum_{i=0}^1\varphi_id\Delta_i^{-1/2}\rho_i+(d\varphi_i)\Delta_i^{-1/2}\rho_i.$$
Let us explain why $\mathfrak{R}$ is a good parametrix for $p<d_S$, i.e. is bounded on $L^p$ if $p<d_S$. First, $d\Delta_0^{-1/2}$ is the Riesz transform on $M_0$, which is bounded by hypothesis. Also, $\varphi_1d\Delta_1^{-1/2}\rho_1$ is a pseudo-differential operator with compact support, and hence is bounded on $L^p$; $(d\varphi_0)\Delta_0^{-1/2}\rho_0$ is a operator with smooth kernel and compact support, hence is bounded on $L^p$. Finally, the operator $(d\varphi_0)\Delta_0^{-1/2}\rho_0$ is bounded on $L^p$ if $p<d_S$, which comes from the facts that $d\varphi_0$ is compactly supported and that for $p<d_S$,

$$\Delta_0^{-1/2} : L^p\rightarrow L^{\frac{np}{n-p}}.$$
The second part of Carron's proof is to show that the error term when we approximate $d\Delta^{-1/2}$ by $\mathfrak{R}$ can be controled on $L^p$ if $p<d_S$. \\\\
In order to improve Carron's result, two things have to be done: first, to find a parametrix for the Riesz transform which is bounded on $L^p$ for $p\geq d_S$, and secondly, to improve the estimates of the error term in order to show that the error term is bounded on $L^p$ for $p\geq d_S$, and not only for $p<d_S$.

\subsection{The case where $M$ has several ends}

In this subsection, we prove Theorem \ref{Riesz p-hyperbolique} in the case where $M$ has several ends. We first remark that the boundedness of the Riesz transform of $M$ on $L^p$ for $p\in [p_0,2]$ is a consequence of Carron's work \cite{Carron4}. We will thus only prove boundedness in the range $[2,\min(d_H,p_1))$. We take the same parametrix for $e^{-t\sqrt{\Delta}}$ than in Carron \cite{Carron4}:

$$\mathfrak{E}(\sigma,u)=\varphi_0e^{-\sigma\sqrt{\Delta_0}}\rho_0 u+\varphi_1e^{-\sigma\sqrt{\Delta_1}}\rho_1 u.$$
The main observation is that when $p\in(p_0,\min(p_1,d_H))$, the corresponding parametrix for the Riesz transform $\mathfrak{R}=d\int_0^\infty \mathfrak{E}(\sigma,\cdot)\,d\sigma$ is bounded on $L^p$. Let us explain this now. We have seen in the previous paragraph that 

$$\mathfrak{R}:=\sum_{i=0}^1\varphi_id\Delta_i^{-1/2}\rho_i+(d\varphi_i)\Delta_i^{-1/2}\rho_i,$$
and that under the hypothesis of Theorem \ref{Riesz p-hyperbolique}, the operators $\varphi_0d\Delta_0^{-1/2}\rho_0$, $\varphi_1d\Delta_1^{-1/2}\rho_1$ and $(d\varphi_1)\Delta_1^{-1/2}\rho_1$ are bounded on $L^p$ for $p\in (p_0,p_1)$. It remains the operator $(d\varphi_0)\Delta_0^{-1/2}\rho_0$. By the fact that $M_0$ satisfies the Sobolev inequality, $M_0$ is $2-$hyperbolic. Thus by the result of Proposition \ref{laplacien non p-parabolique} and interpolation, $(d\varphi_1)\Delta_1^{-1/2}\rho_1$ is bounded on $L^p$ if $p\in [2,d_H)$. Therefore, $\mathfrak{R}$ is bounded for every $p\in(p_0,\min(p_1,d_H))$. All that remains to be done is to show that the corresponding error term is bounded on $L^p$ when $p\in(p_0,\min(p_1,d_H))$, and for this we need to improve the error estimates done in \cite{Carron4}.\\\\
Let $p\in(p_0,\min(p_1,d_H))$; we choose some fixed $q>\frac{d_S}{d_S-2}$ satisfying $p<q<\min(p_1,d_H)$. We will also denote $d:=d_S$. According to \cite{Carron4}, the error term in the parametrix of the Riesz transform is $dg$, where

$$g=\int_0^\infty \int_{\R_+\times R_+\times M}G(\sigma,s,x,y)\left[\left(-\frac{\partial^2}{\partial\sigma^2}+\Delta\right) \mathfrak{E}(\cdot,u)(s,y)\right]d\sigma dsdy,$$ 
$G$ being the Green function of $\left(-\frac{\partial^2}{\partial t^2}+\Delta\right)$ on $M\times \R_+$ with Dirichlet boundary conditions on $M\times\{0\}$.
We let
$$\left(-\frac{\partial^2}{\partial \sigma^2}+\Delta\right)\mathfrak{E}(\sigma,u)=f_0(\sigma,.)+f_1(\sigma,.),$$
where the functions $f_i$ are defined by

$$f_i(t,.)=(\Delta\varphi_i)\left(e^{-t\sqrt{\Delta_i}}\rho_i u\right)-2\langle d\varphi_i,\nabla e^{-t\sqrt{\Delta_i}}\rho_i u\rangle.$$
In \cite{Carron4}, estimates on the $f_i$ are shown. However, since we do not assume $p<d_S$, the corresponding estimates for $f_0$ will not hold in our case. Instead, we will estimate a modified function $\tilde{f}_0$, that we define by

$$\tilde{f}_0(t,.)=\left[\sum_j\mathbf{1}_{A_j}(\Delta\varphi_0)\left(\psi(t)-\left(\psi(t)\right)_{A_j}\right)\right]-2\langle d\varphi_0,\nabla e^{-t\sqrt{\Delta_0}}\rho_0 u\rangle,$$
where $A=\sqcup_j A_j$, each $A_j$ being connected and smooth (see subsection 1 for the definition of $A$), and $\left(\psi(t)\right)_{A_j}$ denotes the average of $\psi$ on $A_j$.
We first show estimates on $f_1$ and $\tilde{f}_0$:

\begin{lem}\label{estimees f}

If $\alpha=d\left(\frac{1}{p}-\frac{1}{q}\right)>0$, then there exists a constant $C$ such that

\begin{equation}\label{estimee f_0}
||\tilde{f}_0(t,.)||_1+||\tilde{f}_0(t,.)||_p\leq \frac{C}{(1+t)^{1+\alpha}}||u||_p,\,\forall t>0.
\end{equation}
and

\begin{equation}\label{estimee f_1}
||f_1(t,.)||_1+||f_1(t,.)||_p\leq \frac{C}{(1+t)^{1+\alpha}}||u||_p,\,\forall t>0.
\end{equation}

\end{lem}
\textit{Proof of Lemma \ref{estimees f}.} We begin by $f_1$. In \cite{Carron4}, it is shown that for some constant $\lambda>0$,

$$||f_1(t,.)||_1+||f_1(t,.)||_p\leq e^{-\lambda t}||u||_p,\,\forall t>0,$$
which of course implies

$$||f_1(t,.)||_1+||f_1(t,.)||_p\leq \frac{C}{(1+t)^{1+\alpha}}||u||_p,\,\forall t>0.$$
Now we turn to $\tilde{f}_0$. Since $d\Delta_0^{-1/2}$ is bounded on $L^q(M_0)$, and $e^{-t\sqrt{\Delta_0}}$ is analytic on $L^r$ for $1<r<\infty$ (see \cite{Stein} or \cite{Coulhon-Saloff-Varopoulos}, this comes from the subordination identity), we have

$$||\nabla e^{-t\sqrt{\Delta_0}}||_{q,q}\leq \frac{C}{t},\,\forall t>0.$$
But

$$||e^{-t\sqrt{\Delta_0}}||_{p,q}\leq \frac{C}{t^{d\left(\frac{1}{p}-\frac{1}{q}\right)}}=\frac{C}{t^\alpha},\,\forall t>0,$$
where $\alpha=d\left(\frac{1}{p}-\frac{1}{q}\right)>0$. We get

$$||\nabla e^{-t\sqrt{\Delta_0}}||_{p,q}\leq ||\nabla e^{-\frac{t}{2}\sqrt{\Delta_0}}||_{q,q}||e^{-\frac{t}{2}\sqrt{\Delta_0}}||_{p,q}\leq \frac{C}{t^{1+\alpha}}.$$
We also have (cf \cite{Carron4})

$$||\nabla e^{-t\sqrt{\Delta_i}}||_{L^p(U)\rightarrow L^q(F)}\leq C,\,\forall t\leq 1,$$
if $U$ is an open subset and $F$ a compact set at positive distance from $U$. Therefore we get

\begin{equation}\label{gradient hors diagonale}
||\nabla e^{-t\sqrt{\Delta_0}}||_{L^p(U)\rightarrow L^q(F)}\leq \frac{C}{(1+t)^{1+\alpha}},\,\forall t>0.
\end{equation}
Since for every $F$ compact, $L^q(F)\hookrightarrow L^1(F)$ and $L^q(F)\hookrightarrow L^p(F)$, and given that the support of $\rho_0$ and $A$ are disjoint, we obtain

$$||\langle d\varphi_0,\nabla e^{-t\sqrt{\Delta_0}}\rho_0 u\rangle||_{L^1}+||\langle d\varphi_0,\nabla e^{-t\sqrt{\Delta_0}}\rho_0 u\rangle||_{L^p}\leq \frac{C}{(1+t)^{1+\alpha}}||u||_p,\,\forall t>0.$$
It remains the term $\left[\sum_j\mathbf{1}_{A_j}(\Delta\varphi_0)\left(\psi(t)-\left(\psi(t)\right)_{A_j}\right)\right]$. We have, by the Poincar\'{e} inequality $L^q$ on $A_j$:

$$||\psi(t)-\left(\psi(t)\right)_{A_j}||_{L^q(A_j)}\leq C||\nabla e^{-t\sqrt{\Delta_0}}\rho_0 u||_{L^q(A_j)}\leq \frac{C}{(1+t)^{1+\alpha}}||u||_p,\,\forall t>0.$$
Hence the estimates for $\tilde{f}_0$.

\cqfd
Now, we decompose $g$ into $g_1+g_2$, with

$$g_1(x)=\int_{\R_+\times\R_+\times M}G(\sigma,s,x,y)\tilde{f}_0(s,y)d\sigma dsdy+\int_{\R_+\times\R_+\times M}G(\sigma,s,x,y)f_1(s,y)d\sigma dsdy,$$
and

$$g_2(x)=\sum_j\int_{\R_+\times\R_+\times M}G(\sigma,s,x,y)\mathbf{1}_{A_j}(y)(\Delta\varphi_0)(y)\left(\psi(s)\right)_{A_j}))d\sigma dsdy.$$
We have, in an equivalent way,

$$g_1=\frac{2}{\sqrt{\pi}}\int_{\R_+\times\R_+}e^{-r^2}\left(\int_0^{\frac{s^2}{4r^2}}(e^{-t\Delta}(\tilde{f}_0(s,.)+f_1(s,.))dt\right)drds,$$
and

$$g_2=\sum_j\frac{2}{\sqrt{\pi}}\int_{\R_+\times \R_+}e^{-r^2}\left(\int_0^{\frac{s^2}{4r^2}}e^{-t\Delta}(\mathbf{1}_{A_j}(\Delta\varphi_0)\left(\psi(s)\right)_{A_j})dt\right)drds.$$
In order to conclude the proof of Theorem \ref{Riesz p-hyperbolique} in the case of several ends, we have to show that $||dg_1||_p+||dg_2||_p\leq C||u||_p.$ This will be done in the next two Lemmas. Let us begin by

\begin{lem}\label{erreur1}

There exists a constant $C$ such that for every $u\in L^p$,

$$||dg_1||_p\leq C||u||_p.$$

\end{lem}
\pr According to Proposition 2.1 in \cite{Carron4}, it is enough to show that $||g_1||_p+||\Delta g_1||_p\leq C||u||_p$. The term $||\Delta g_1||_p$ is the easiest: defining $h:=\tilde{f}_0+f_1$, we have

$$\begin{array}{rcl}
\Delta g_1&=&\frac{2}{\sqrt{\pi}}\int_{\R_+\times\R_+} e^{-r^2}\left(\int_0^{\frac{s^2}{4r^2}}\Delta(e^{-t\Delta}h(s,.))dt\right) drds\\\\
	&=&-\frac{2}{\sqrt{\pi}}\int_{\R_+\times\R_+} e^{-r^2}\left(\int_0^{\frac{s^2}{4r^2}}\frac{d}{dt}(e^{-t\Delta}h(s,.))dt\right) drds\\\\
	&=&\frac{2}{\sqrt{\pi}}\int_{\R_+\times\R_+} e^{-r^2}\left(h(s,.)-(e^{-\frac{s^2}{4r^2}\Delta}h(s,.))\right)drds.
\end{array}$$
Hence, by \eqref{estimee f_0} and \eqref{estimee f_1},

$$\begin{array}{rcl}
||\Delta g_1||_p&\leq& \frac{4}{\sqrt{\pi}}\int_{\R_+\times\R_+}e^{-r^2}||h(s,.)||_pdrds\\\\
	&\leq& \frac{4}{\sqrt{\pi}}\left(\int_{\R_+\times\R_+}e^{-r^2}\frac{C}{(1+s)^{1+\alpha}}drds\right) ||u||_p\\\\
	&\leq& C||u||_p
\end{array}$$
For $||g_1||_p$, using

$$||e^{-t\Delta}||_{1,p}\leq \frac{C}{t^{\frac{d}{2}\left(1-\frac{1}{p}\right)}},$$
and \eqref{estimee f_0}, \eqref{estimee f_1}, we have

$$\begin{array}{rcl}
||g_1||_p&\leq& \frac{2}{\sqrt{\pi}}\int_{\R_+\times\R_+}e^{-r^2}\left(\int_0^{\frac{s^2}{4r^2}}||e^{-t\Delta}h(s,.)||_pdt\right)dsdr\\\\
	&\leq& \frac{2}{\sqrt{\pi}}\left(\int_{\R_+\times\R_+}e^{-r^2}\left(\int_0^{\frac{s^2}{4r^2}}\frac{C}{\max\left(1,t^{\frac{d}{2}(1-\frac{1}{p})}\right)(1+s)^{1+\alpha}}dt\right)dsdr\right)||u||_p\\\\
	&\leq& C\left(\int_{\R_+\times\R_+}e^{-r^2}\frac{1}{\max\left(1,t^{\frac{d}{2}(1-\frac{1}{p})}\right)(1+2r\sqrt{t})^{\alpha}}dtdr\right)||u||_p
\end{array}$$
We seperate the integral in $\int_{t\leq r^{-2}}+\int_{t\geq r^{-2}}=I_1+I_2$. The integral $I_1$ is finite if and only if

$$(-2)\left(\frac{d}{2}\left(1-\frac{1}{p}\right)-1\right)<1,$$
which is equivalent to

$$p>\frac{d}{d-1}.$$
Since $p>p_0>\frac{d}{d-1}$, this is automatically satisfied. For $I_2$,

$$\begin{array}{rcl}
I_2&\leq&\int_0^\infty e^{-r^2}\left(\int_{r^{-2}}^\infty \frac{1}{t^{\frac{d}{2}(1-\frac{1}{p})}}\frac{1}{(r\sqrt{t})^\alpha}dt\right)dr\\\\
	&\leq& \int_0^\infty e^{-r^2}\frac{1}{r^\alpha}\left(\int_{r^{-2}}^\infty\frac{1}{t^{\frac{d}{2}(1-\frac{1}{p})}}\frac{1}{(\sqrt{t})^\alpha}dt\right)dr
\end{array}$$
The integral in $t$ is finite if and only if

$$\frac{d}{2}\left(1-\frac{1}{p}\right)+\frac{\alpha}{2}>1,$$
and recalling that $\alpha=d\left(\frac{1}{p}-\frac{1}{q}\right)$, we find that it is equivalent to

$$q>\frac{d}{d-2},$$
which is satisfied by assumption. The integral in $r$ is then

$$\int_0^\infty e^{-r^{2}}\frac{1}{r^{\alpha-2\left(\frac{d}{2}(1-\frac{1}{p})+\frac{\alpha}{2}-1\right)}}dr,$$
which is finite if and only if

$$\alpha-d\left(1-\frac{1}{p}\right)-\alpha+2<1,$$
which is equivalent to

$$p>\frac{d}{d-1},$$
which is satisfied by assumption since $p>p_0>\frac{d}{d-1}$.
\cqfd
Now we turn to estimate $g_2$, which will conclude the proof of Theorem \ref{Riesz p-hyperbolique} in the case of several ends.

\begin{lem}

$$||dg_2||_p\leq C||u||_p.$$

\end{lem}
\pr According to Proposition 2.1 in \cite{Carron4}, it is enough to show that $||g_2||_p+||\Delta g_2||_p\leq C||u||_p$. We begin to show that $||g_2||_p\leq C||u||_p$. We have

$$\begin{array}{rcl}
g_2(x)&=&\sum_j\frac{2}{\sqrt{\pi}}\int_{\R_+\times\R_+}e^{-r^2}\left(\int_0^{\frac{s^2}{4r^2}}e^{-t\Delta}\left(\mathbf{1}_{A_j}(\Delta\varphi_0)\left(\psi(s)\right)_{A_j}\right)(x)dt\right)drds\\\\
	&=&\sum_j\frac{2}{\sqrt{\pi}}\int_{\R_+\times\R_+}e^{-r^2}\left(\int_0^{\frac{s^2}{4r^2}}\left(\psi(s)\right)_{A_j}e^{-t\Delta}\left(\mathbf{1}_{A_j}\Delta\varphi_0\right)(x)dt\right)drds,
\end{array}$$
therefore

$$||g_2||_p\leq \sum_j\frac{2}{\sqrt{\pi}}\int_{\R_+\times\R_+}e^{-r^2}\left(\int_0^{\frac{s^2}{4r^2}}\left(\frac{1}{|A_j|}\int_{A_j}e^{-s\sqrt{\Delta_0}}|\rho_0 u|\right)||e^{-t\Delta}\chi||_pdt\right)drds,$$
where we have defined $\chi:=\Delta\varphi_0=\Delta(\varphi_0-1)$. Using the fact that $||e^{-t\Delta}||_{1,p}\leq \frac{C}{t^{\frac{d}{2}\left(1-\frac{1}{p}\right)}}$, the analyticity of $e^{-t\Delta}$ on $L^p$, and the fact that $\varphi_0-1$ is smooth with compact support,

$$||e^{-t\Delta}\chi||_p\leq \frac{C}{\max\left(1,t^{1+\frac{d}{2}\left(1-\frac{1}{p}\right)}\right)},\,\forall t>0.$$
Furthermore, we have for every  $p>1$,

$$1+\frac{d}{2}\left(1-\frac{1}{p}\right)>1,$$
and consequently

$$\int_0^\infty ||e^{-t\Delta}\chi||_pdt<\infty.$$
So

$$\begin{array}{rcl}
||g_2||_p&\leq& C\sum_j\int_{\R_+\times\R_+}e^{-r^2}\left(\int_{A_j}e^{-s\sqrt{\Delta_0}}|\rho_0 u|\right)drds\\\\
	&\leq& C\sum_j\int_{A_j}\left(\int_0^\infty e^{-s\sqrt{\Delta_0}}|\rho_0 u|ds\right) \\\\
	&\leq& C\sum_j\int_{A_j} \Delta_0^{-1/2}|\rho_0 u|
\end{array}$$
According to Proposition (\ref{laplacien non p-parabolique}), $\Delta_0^{1/2} : L^p\rightarrow L^p_{loc}\hookrightarrow L^1_{loc}$, which implies that

$$||g_2||_p\leq C||u||_p.$$
Let us turn now to $\Delta g_2$: as for $g_1$, we have

$$\begin{array}{rcl}
||g_2||_p&\leq& \sum_j\frac{4}{\pi}\int_{\R^2_+}e^{-r^2}||\left(\psi(s)\right)_{A_j}(\Delta\varphi_0)||_pdrds\\\\
	&\leq& C\sum_j\int_0^\infty|\left(\psi(s)\right)_{A_j}|ds,
\end{array}$$
and by the argument we have already used,

$$\sum_j\int_0^\infty|\left(\psi(s)\right)_{A_j}|ds\leq C||u||_p,$$
which concludes the proof.
\cqfd

\subsection{The case where $M$ has one end}

In this subsection, we prove Theorem \ref{Riesz p-hyperbolique} when $M$ has \textit{only one end}. As we already explained, the parametrix $\mathfrak{R}$ for the Riesz transform constructed by Carron has a term which is not bounded on $L^p$ when $p>d_H$: more precisely, the term $(d\varphi)\Delta_0^{-1/2}\rho_0$ is not bounded on $L^p$ if $p>d_H$. Hence, we have to modify the parametrix. The main idea is the following: notice that since $M$ has only one end, $d\varphi$ is the supported in $A$ which is a \textit{connected} annulus. Thanks to the $L^p$ Poincar\'{e} inequality in $A$, there is a constant $C$ such that

$$\left|\left|v-\frac{1}{|A|}\int_{A}v\right|\right|_{L^p(A)}\leq C||\nabla v||_p\,\forall v\in C^\infty(A).$$
Applying this to $\Delta_0^{-1/2}\rho_0 u$ for $u\in C_0^\infty(M)$, we get for $p\in (p_0,p_1)$

$$\left|\left| \Delta_0^{-1/2}\rho_0 u-\frac{1}{|A|}\int_A \Delta_0^{1/2}\rho_0 u\right|\right|_{L^p(A)}\leq C||\nabla \Delta_0^{-1/2}\rho_0 u||_p\leq C||u||_p,$$
where in the last inequality we used the fact that the Riesz transform on $M_0$ is bounded on $L^p$ if $p\in (p_0,p_1)$. This implies that the modified parametrix

$$\mathfrak{S}u=\sum_{i=0}^1\varphi_iT_i\rho_i u+(d\varphi_1)\Delta_1^{-1/2}\rho_1u+(d\varphi_0)\left(\Delta^{-1/2}\rho_0 u-\left(\frac{1}{|A|}\int_{A}\Delta_0^{-1/2}\rho_0u\right)\right),$$
where $T_i:=d\Delta_i$, is bounded on $L^p$ for every $p\in (p_0,p_1)$. The corresponding parametrix for $e^{-t\sqrt{\Delta}}$ is given by

$$\mathfrak{F}(\sigma,u)=\mathfrak{E}(\sigma,u)-(\varphi_0-1)\left(\frac{1}{|A|}\int_{A}e^{-\sigma \sqrt{\Delta_0}}\rho_0u\right),$$
i.e.

$$\mathfrak{S}u=d\int_0^\infty \mathfrak{F}(\sigma,u) \,d\sigma.$$
The supplementary term that we have added to the parametrix of $e^{-t\sqrt{\Delta}}$ is 

$$-(\varphi_0-1)\left(\frac{1}{|A|}\int_{A}e^{-\sigma \sqrt{\Delta_0}}\rho_0u\right),$$
which vanishes when $\sigma=0$, since $A$ and the support of $\rho_0$ are disjoint by hypothesis. So we have, as should be,

$$\mathfrak{F}(0,u)=u.$$
Notice also that since $\varphi_0-1$ is compactly supported, the integral with respect to $\sigma$ of this supplementary term is analogous to the term $G_3$ in the parametrix of $\Delta^{-1/2}$ constructed By Carron-Coulhon-Hassell in \cite{Carron-Coulhon-Hassell}: its kernel $k(x,y)$ is non-zero only if $x$ is in $K_3$ and $y$ is in $M\setminus K_1$.\\

We thus have constructed a parametrix $\mathfrak{S}$ for the Riesz transform, which is bounded on $L^p$ for $p\in (p_0,p_1)$. As in the proof of Theorem \ref{Riesz p-hyperbolique} in the case where $M$ has several ends, it remains to show that the error term is also bounded on $L^p$.\\\\
We will use the calculations made in the previous subsection. This time, we have (with $f_1$ defined as in the previous subsection)

$$\begin{array}{rcl}\left(-\frac{\partial^2}{\partial \sigma^2}+\Delta\right)\mathfrak{F}(\sigma,u)&=&f_1(\sigma,.)+(\Delta\varphi_0)\left(e^{-t\sqrt{\Delta}}\rho_0 u\right)-2\langle d\varphi_0,\nabla e^{-t\sqrt{\Delta_0}}\rho_0 u\rangle\\\\
&&-(\Delta\varphi_0)\left(\frac{1}{|A|}\int_{A}e^{-\sigma\sqrt{\Delta_0}}\rho_0 u\right)-(\varphi_0-1)\left(\frac{1}{|A|}\int_{A}\Delta_0 e^{-\sigma\sqrt{\Delta_0}}\rho_0 u\right).
\end{array}$$
Define as in the previous section

$$\tilde{f}_0(\sigma,.):=(\Delta\varphi_0)\left(e^{-t\sqrt{\Delta}}\rho_0 u-\left(\frac{1}{|A|}\int_{A}e^{-\sigma\sqrt{\Delta_0}}\rho_0 u\right)\right)-2\langle d\varphi_0,\nabla e^{-t\sqrt{\Delta_0}}\rho_0 u\rangle,$$
and also

$$\bar{f}_0(\sigma,.)=(\varphi_0-1)\left(\frac{1}{|A|}\int_{A}\Delta_0 e^{-\sigma\sqrt{\Delta_0}}\rho_0 u\right).$$
We have the following estimates on $f_1$, $\tilde{f}_0$ and $\bar{f}_0$:

\begin{lem}\label{estimee f 2}
If $\alpha=d\left(\frac{1}{p}-\frac{1}{q}\right)$, then for all $\sigma>0$,

$$||f_1(\sigma,\cdot)||_1+||f_1(\sigma,\cdot)||_p\leq \frac{C}{(1+\sigma)^{1+\alpha}}||u||_p,$$

$$||\tilde{f}_0(\sigma,\cdot)||_1+||\tilde{f}_0(\sigma,\cdot)||_p\leq \frac{C}{(1+\sigma)^{1+\alpha}}||u||_p,$$
and
$$||\bar{f}_0(\sigma,\cdot)||_1+||\bar{f}_0(\sigma,\cdot)||_p\leq \frac{C}{(1+\sigma)^{1+\alpha}}||u||_p.$$

\end{lem}
Once this Lemma is established, the estimate of the error term proceeds as in the proof of Theorem \ref{Riesz p-hyperbolique} in the case where $M$ has more than one end. All we have to do is thus to prove the above estimates.\\\\
\textit{Proof of Lemma \ref{estimee f 2}:} We already proved the estimates on $f_1$ and $\tilde{f}_0$ in Lemma \ref{estimees f}.
It remains to treat $\bar{f}_0$. First, by analyticity of $e^{-t\sqrt{\Delta_0}}$,

\begin{equation}\label{terme d'erreur analytique}
\left|\left|\Delta_0e^{-t\sqrt{\Delta_0}}\right|\right|_{p,p}\leq \frac{C}{t^2},
\end{equation}
and therefore, using the fact that $\bar{f}_0(\sigma,\cdot)$ has compact support independant of $u$,

$$||\bar{f}_0(\sigma,\cdot)||_{1,1}+||\bar{f}_0(\sigma,\cdot)||_{p,p}\leq \frac{C}{\sigma^{2}}.$$
The proof will be complete once we show that $\Delta_0 e^{-t\sqrt{\Delta_0}}$ is a bounded operator $L^p(M_0\setminus A_\delta)\rightarrow L^\infty(A)$ when $t\rightarrow 0$ (where $\delta$ is a strictly positive constant, and where $A_\delta$ is the set of points whose distance to $A$ is greater than $\delta$). For this, we use the subordination identity:

$$e^{-\sigma \sqrt{\Delta_0}}=\frac{\sigma}{2\sqrt{\pi}}\int_0^\infty e^{-\frac{\sigma^2}{4t}}e^{-t\Delta_0}\frac{dt}{t^{3/2}},$$
so that

\begin{equation}\label{terme d'erreur un bout}
\Delta_0 e^{-\sigma \sqrt{\Delta_0}}=-\frac{\sigma}{2\sqrt{\pi}}\int_0^\infty e^{-\frac{\sigma^2}{4t}}\left(\frac{\partial}{\partial t}e^{-t\Delta_0}\right)\frac{dt}{t^{3/2}}.
\end{equation}
According to \cite{Davies2}, Corollary 5 (se also \cite{Saloff-Coste1}, Theorem 5.2.15), the Sobolev inequality on $M_0$ implies

\begin{equation}\label{estimee derivee temps 1}
\left|\frac{\partial p^0_t(x,y)}{\partial t}\right|\leq \frac{C}{t^{n+1}}e^{-c\frac{d^2(x,y)}{t}},\,\forall (x,y)\in M_0\times M_0,\,\forall t>0,
\end{equation}
where $p^0_\sigma(x,y)$ is the heat kernel on $M_0$. So, if $\Omega$ is an open set and $F$ a compact set such that $d(F,\Omega)\geq\varepsilon>0$, then,

\begin{equation}\label{estimee derivee temps 2}
\left|\frac{\partial p^0_t(x,y)}{\partial t}\right|\leq \frac{C}{t^{n+1}}\exp\left(-c\frac{\varepsilon^2}{t}\right),\,\forall t>0,\,\forall x\in F,\,\forall y\in\Omega.
\end{equation}
The estimates \eqref{estimee derivee temps 1} and \eqref{estimee derivee temps 2} imply the existence of a constant (depending of the lower bound on the Ricci curvature of $M$ and of $\delta$) such that, if $t\leq 1$,

\begin{equation}\label{estimee Lp L infty}
\left|\left|\frac{\partial}{\partial t}e^{-t\Delta_0}\right|\right|_{L^p(M_0\setminus A_\delta)\rightarrow L^\infty(A)}\leq C.
\end{equation}
Indeed, denoting $k_t(x,y)=\frac{1}{t^{n+1}}\exp\left(-c\frac{d^2(x,y)}{t}\right)$, and $K_t$ the operator with kernel $k_t$, then

\begin{equation}\label{estimee L1 Linfty}
K_t : L^1(\Omega)\rightarrow L^\infty(F)
\end{equation}
is uniformly bounded when $t\rightarrow 0$: this comes from the fact that by \eqref{estimee derivee temps 2},

$$||K_t||_{ L^1(\Omega)\rightarrow L^\infty(F)}=\sup_{x\in F,\,y\in \Omega}k_t(x,y)\leq C,\,\forall t\leq 1.$$
Furthermore,

\begin{equation}\label{estimee Linfty Linfty}
K_t : L^\infty(\Omega)\rightarrow L^\infty(F)
\end{equation}
is uniformly bounded when $t\rightarrow 0$. To show this, we have to prove that

$$\sup_{x\in F}\int_\Omega k_t(x,y)\leq C,\,\hbox{for all }t\hbox{ small enough.}$$
But for $t\leq 1$ and $x\in F$, $y\in\Omega$, \eqref{estimee derivee temps 1} yields

\begin{equation}\label{estimee derivee temps 3}
k_t(x,y)\leq C_1\exp\left(-\frac{c}{2}\frac{d^2(x,y)}{t}\right).
\end{equation}
We then use the fact that the volume of balls of radius $r$ is bounded by $e^{ar}$ for a certain constant $a$, since the Ricci curvature is bounded from below on $M$; therefore, we deduce that if $t$ is small enough so that for every $x\in F$, $y\in \Omega$,

$$\frac{c}{2}\frac{\varepsilon}{t}> a,$$
then by \eqref{estimee derivee temps 3},

$$\sup_{x\in F}\int_\Omega k_t(x,y)\leq C_2.$$
Finally, \eqref{estimee Lp L infty} is obtained by interpolation from \eqref{estimee L1 Linfty} and from \eqref{estimee Linfty Linfty}. Using in addition the analyticity of $e^{-t\Delta_0}$ and the fact that $e^{-\frac{1}{2}\Delta_0} : L^p\rightarrow L^\infty$, we obtain that

$$\left|\left|\frac{\partial }{\partial t}e^{-t\Delta_0}\right|\right|_{L^p(M_0\setminus A_\delta)\rightarrow L^\infty(A)}=\left|\left|\Delta_0 e^{-t\Delta_0}\right|\right|_{L^p(M_0\setminus A_\delta)\rightarrow L^\infty(A)}\leq \frac{C}{1+t},\,\forall t>0.$$
In particular,

$$\left|\left|\frac{\partial }{\partial t}e^{-t\Delta_0}\right|\right|_{L^p(M_0\setminus A_\delta)\rightarrow L^\infty(A)}\leq C,\,\forall t>0.$$
Using \eqref{terme d'erreur un bout}, we then obtain
$$\left|\left| \Delta_0 e^{-\sigma\sqrt{\Delta_0}}\right|\right|_{L^p(M_0\setminus A_\delta)\rightarrow L^\infty(A)}\leq C,\,\forall \sigma>0,$$
and reminding of \eqref{terme d'erreur analytique}, we have

$$\left|\left| \Delta _0e^{-\sigma\sqrt{\Delta_0}}\right|\right|_{L^p(M_0\setminus A_\delta)\rightarrow L^p(A)}+\left|\left| \Delta _0e^{-\sigma\sqrt{\Delta_0}}\right|\right|_{L^p(M_0\setminus A_\delta)\rightarrow L^1(A)}\leq \frac{C}{(1+\sigma)^2},\,\forall \sigma>0.$$
This implies that $\bar{f}_0(\sigma,\cdot)$ is bounded as an operator on $L^p$ when $\sigma\to0$. Using the fact that the support of $\bar{f}_0(\sigma,\cdot)$ is compact and independant of $u$, we get

$$||\bar{f}_0(\sigma,\cdot)||_{1,1}+||\bar{f}_0(\sigma,\cdot)||_{p,p}\leq \frac{C}{(1+\sigma)^{1+\alpha}}.$$

\cqfd

\subsection{Proof of the corollaries to Theorem \ref{Riesz p-hyperbolique}}
In this final subsection, we give the proofs of Corollaries \ref{asymptotic conic}, \ref{Riesz Ricci positive hors d'un compact}, \ref{Lie group} and \ref{Riesz non borne somme connexe p-hyperbolique}.\\\\
\textit{Proof of Corollary \ref{asymptotic conic}:} Using the result of H.Q. Li \cite{Li1} and noticing that the conic manifold $M_0$ satisfies $d_S=$dim$(M_0)>2$ and that $p_0>d_S$, we can apply Theorem \ref{Riesz p-hyperbolique} to get that the Riesz transform on $M$ is bounded on $2\leq p<p_0$. The boundedness on $L^p$ of the Riesz transform on $M$ for $1<p<2$ follows from Coulhon-Duong's result \cite{Coulhon-Duong2} and the fact that $M$ satisfies a Sobolev inequality. Now, if the Riesz transform on $M$ were bounded on $L^p$ for $p\in (1,p_1)$ with $p_1$, then applying Theorem \ref{Riesz p-hyperbolique} reversing the roles of $M_0$ and $M$, we would get that the Riesz transform on $M_0$ is bounded on $L^p$ for every $p\in (1,p_1)$, which is false by H.Q. Li's result.
\cqfd
\textit{Proof of the Corollary \ref{Riesz Ricci positive hors d'un compact}:} The hypothesis on the volume of balls implies (see \cite{Carron4}) that for every compact set $K$ in $M$, for all $1\leq p\leq q\leq \infty$, and for every $t\geq1$,

$$\left|\left| e^{-t\Delta}\right|\right|_{L^p(K)\rightarrow L^q(M)}\leq \frac{C_K}{t^{\nu\left(\frac{1}{p}-\frac{1}{q}\right)}},$$
and the proof of Theorem (\ref{Riesz p-hyperbolique}) applies.
\cqfd
\textit{Proof of the Corollary \ref{Lie group}:}  It is known by \cite{Alexopoulos} that the Riesz transform on $M_0$ is bounded on $L^p$ for every $1<p<\infty$. The Sobolev inequality on a simply connected, nilpotent Lie group is proved in \cite{Coulhon-Saloff-Varopoulos}, p.56. The boundedness on $L^p$ of the Riesz transform on $M$ for $1<p<2$ follows from Coulhon-Duong's result \cite{Coulhon-Duong2} and the fact that $M$ satisfies a Sobolev inequality. Furthermore, the simple connectedness of $M_0$ implies that it has only one end. Therefore, we can apply Theorem (\ref{Riesz p-hyperbolique}) to get the result for $2\leq p<\infty$.
\cqfd
\textit{Proof of Corollary \ref{Riesz non borne somme connexe p-hyperbolique}:} By interpolation, it is enough to prove that the Riesz transform on $M$ is not bounded on $L^p$ for $n<p<q$. We proceed by contradiction: let us assume that the Riesz transform on $M$ is bounded on $L^p$ for a certain $n<p<q$. Then, since $M$ is $q-$hyperbolic according to Corollary (\ref{somme connexe p-hyperbolique}), by applying Theorem (\ref{Riesz p-hyperbolique}) we find that the Riesz transform on $M\#M$ is bounded on $L^r$, for some $n<r<p$. But $M\#M=(\R^n\#\R^n)\#(N\#N)$, and since $M\#M$ is also $q-$hyperbolic, Theorem (\ref{Riesz p-hyperbolique}) implies that the Riesz transform on the disjoint union of $\R^n\#\R^n$ and of $N\#N$ is bounded on $L^s$, for some $n<s<r$. But we know, according to \cite{Carron-Coulhon-Hassell} that the Riesz transform on $\R^n\#\R^n$ is not bounded on $L^s$ if $s\geq n$; hence a contradiction.
\cqfd

\begin{center}{\bf Acknowledgments} \end{center}
This article is part of the PhD thesis of the author. The author would like to thank his advisor G. Carron, for inspiring discussions and support.

\bibliographystyle{plain}

\bibliography{bibliographie}

\end{document}